\documentclass[psfig]{article}

\usepackage{amssymb, amsmath, amsthm}

\newtheorem{theo}{Theorem}[section]
\newtheorem{lem}{Lemma}[section]

\newcommand{\BB}[1]{\ensuremath{\mathbb{#1}}}
\newcommand{\N}{\ensuremath{\BB{N}}}

\newcommand{\R}{\ensuremath{\BB{R}}}
\newcommand{\Z}{\ensuremath{\BB{Z}}}

\newcommand{\be}{\begin{equation}}
\newcommand{\ee}{\end{equation}}
\newcommand{\bes}{\begin{equation*}}
\newcommand{\ees}{\end{equation*}}
\newcommand{\bi}{\begin{itemize}}
\newcommand{\ei}{\end{itemize}}
\newcommand{\bea}{\begin{eqnarray}}
\newcommand{\eea}{\end{eqnarray}}
\newcommand{\beas}{\begin{eqnarray*}}
\newcommand{\eeas}{\end{eqnarray*}}

\begin{document}
\title{Multiplicity and regularity of periodic solutions for a class of degenerate semilinear wave equations.
}
\author{Jean-Marcel Fokam \\
{\small \tt fokam\,@aun.edu.ng}\\
School of Arts and Sciences, American University of Nigeria\\
Yola,
}


\date{} 



\maketitle

\begin{abstract}
We prove the existence of infinitely many classical periodic solutions for a class of degenerate semilinear wave equations:
\[
u_{tt}-u_{xx}+|u|^{s-1}u=f(x,t),
\]

for all $s>1$. In particular we prove the existence of infinitely many classical solutions for the case $s=3$ posed by Br\'ezis in \cite{BrezisBAMS}.
The proof relies on a new upper a priori estimate, for minimax values of, a perturbed from symmetry, strongly indefinite functional depending on a small parameter.

\end{abstract}

\newpage

\section{Introduction}

\noindent In this paper we construct infinitely many classical time-periodic solutions
for the following semilinear degenerate wave equation with time-dependent forcing
term $f$:
\begin{equation}
    u_{tt}-u_{xx}+g(u)-f(t,x)=0 \label{papier}
\end{equation}
\be
    u(0,t)=u(\pi,t)=0.
\label{bdry}
\ee
where $g(u)=|u|^{s-1}u$ and $F(x,t,u)=g(u)-f(x,t)$, where $f$ is of class $C^2$ and satisfies the Dirichlet boundary conditions.

{\bf Br\'ezis problem}\cite{BrezisBAMS}:\textit{It seems reasonable to conjecture that when $g(u)=u^3$ problem (\ref{papier}),(\ref{bdry}) possesses a solution -even infinitely many solutions- for every $f$(or at least a dense set of $f$'s.)}

\begin{theo}{\rm{ If $f\in C^2$ then there exists infinitely many}} classical solutions of (\ref{papier}),(\ref{bdry}) {\rm{for all}} $s>1$.
\label{conjecture}
\end{theo}
Theorem \ref{conjecture} also prove the existence of classical solutions for a question of Bahri-Berestycki in \cite{BahriBerestyckiActa} on the existence of infinitely many solutions of (\ref{papier}),(\ref{bdry}) for the class of odd nonlinearities, $g(u)=|u|^{s-1}u$.

The \textit{weak version} of the conjecture of Br\'ezis, the existence of \textit{weak solutions for a dense set of} $f$'s has been shown to be true by Tanaka in \cite{Tanaka86}. The problem (\ref{papier}),(\ref{bdry}), \textit{for a given} $f$, has been studied by Tanaka \cite{Tanaka88}, Bartsch-Ding-Lee \cite{BDL}, for arbitrary $s>1$, and Bolle-Ghoussoub-Tehrani \cite{BGT}, Ollivry \cite{Oll} for the case $1<s<2$ however only \textit{weak solutions} have been obtained. As already noticed in \cite{R71} there are two classes of monotone functions for problem (\ref{papier}),(\ref{bdry}), the strongly monotone $F$, $\frac{\partial F}{\partial u}\geq \alpha>0$ which can be compared to the uniformly elliptic case and the degenerate monotone case which allows $\frac{\partial F}{\partial u}=0$. These two classes of monotone functions have been extensively studied by Torelli\cite{Torelli},Rabinowitz \cite{R71}, Hall\cite{Hall}, Hale\cite{Hale}, in the small perturbative case, i.e. with a smallness assumption on $f$ . No such a smallness assumption is assumed here and the result we prove is a global one.

The difficulty in proving the regularity of the weak solutions obtained by \cite{Tanaka88},\cite{BDL},\cite{BGT} lies in the strong monotonicity assumption which is required by the regularity approach of Br\'ezis-Nirenberg, \cite{BN2}. In \cite{BN2} Br\'ezis and Nirenberg show that an $L^\infty$ weak solution is smooth as long as $F$ is smooth and satisfies the strong monotonicity assumption $\frac{\partial F}{\partial u}\geq \varepsilon>0$ which fails here as $g(u)$ has a vanishing derivative.  Note that in the highly degenerate case where $F$ vanishes in an interval, weak solutions in $L^\infty$ need not to be smooth, see \cite{BN2} or \cite{BN} theorem $I.8$. Therefore, to find classical periodic solutions we will proceed differently. In \cite{R78} Rabinowitz developed a regularity theory for this type of  degeneracy where $\frac{\partial F}{\partial u}=0$ is allowed but $g$ strictly monotone ($z_1>z_2$ implies $g(z_1)>g(z_2)$) for equations of the type  (\ref{papier}),(\ref{bdry}) and with $f=0$. The approach in \cite{R78} consisted in seeking viscous approximative solutions, studying a modified equation analogue of (\ref{mpapier}) with $f=0$:
\be
   w_{tt}(\beta)-w_{xx}(\beta)=-|u|^{s-1}u(\beta)+\beta v_{tt}(\beta)
\label{mpapier78}
\ee
(Here $u(\beta)=v(\beta)+w(\beta)$ and $v(\beta)$ is the component of $u(\beta)$ in the direction of the infinite dimensional kernel of $\Box$, with the Dirichlet-periodic boundary conditions. The solution $u$ is split in such a way to tackle the problem stemming from the infinite dimensional kernel of $\Box$.)
with the {\textit{small parameter}} $\beta$ and obtaining compactness via upper priori estimates independently of $\beta$ of the critical values of the modified problem (\ref{mpapier78}), enabling him to send $\beta$ to $0$ and then finding classical solutions. However the problem here contains the forcing term $f$ and the natural functional associated with the problem (\ref{papier}) is no longer even thus  the minimax sets for finding critical values in \cite{R78} do not apply for forced vibrations.

In the eighties and nineties a perturbation theory for this type of problems \textit{-perturbation from symmetry-} was developed,by Bahri-Berestycki \cite{BahriBerestycki81},Bahri-Lions \cite{BahriLions88},Tanaka \cite{Tanaka89}  Struwe \cite{Struwe}, Rabinowitz \cite{R82} and Bolle \cite{Bolle}. The approaches consist in finding growth estimates on some minimax values,$b_n$, and if they grow fast enough, will imply the existence of critical values of the perturbed functional.
Hence it is therefore natural to try to implement these approaches, to tackle the regularity issues stemming from the degenerate monotone semilinear term $g(u)$ and the infinite dimensional kernel of $\Box$ under Dirichlet boundary conditions, to the modified equation, seeking viscous approximative solutions:
\be
   w_{tt}(\beta)-w_{xx}(\beta)=-|u|^{s-1}u(\beta)+\beta v_{tt}(\beta)+f(t,x).
\label{mpapier}
\ee
However the approaches by \cite{BahriBerestycki81},\cite{BahriLions88},\cite{Bolle},\cite{Struwe},\cite{R82},\textit{do not provide an upper explicit upper estimates} on the critical values, and this lead to serious difficulties to obtain compactness of $u(\beta)$, as $\beta\rightarrow 0$.

For even functionals, the identity map is an admissible function in the set of maps considered for the minimax procedure. Information gleaned from the identity map in \cite{R84} has lead to explicit a priori estimates and hence compactness for free vibrations. For forced vibration such an explicit map is lacking and to overcome these difficulties {\textit{we construct a map in the minimax sets of Rabinowitz \cite{R82}, whose energy in $J_{\beta}$ is controlled independently of $\beta$}}. The additional estimate thus obtained lead to the needed compactness needed to pass to the limit as $\beta\rightarrow 0$.

Having constructed minimax values $c_n^m(\delta)$ with upper a priori estimates independently of, the Galerkin parameter $m$ and $\beta$, we need information on the growth of some minimax values $b_n^m$ to show that the $c_n^m(\delta)$ are critical values. To obtain the lower estimates of the growth of the $b_n^m$ we employ the functional $K$ introduced by Tanaka in \cite{Tanaka88} and the Borsuk-Ulam lemma of Tanaka \cite{Tanaka88}, see lemma \ref{Borsuk-Ulam-Tanaka}.

Another advantage of our approach is that it simplifies the weak solutions approach of \cite{Tanaka88}.  In \cite{Tanaka88} some technical lemmas are employed to get information on the index of the weak solution $u$, obtained by passing to the limit in the Galerkin parameter $m$, the index of the critical value of the approximate solution $u_n^m$, obtained from the Galerkin scheme. Here the upper estimate on $c_n^m(\delta)$ is also independent of $m$ thus it allows to simplify the passage to the limit as $m\rightarrow\infty$.

Once the compactness of the sequence $u(\beta)$ is obtained, the regularity will follow by the adapting the argument of \cite{R78} to the problem considered here, in presence of a forcing term $f(x,t)$.\\

Remark:\textit{Upper estimates for criticial values via the approach of \cite{Bolle} and under Dirichlet boundary conditions are in \cite{CDHL} by Castro,Ding and Hernandez-Linares, and  Castro and Clapp \cite{CastroClapp}, for perturbation of a differential operator,the Laplacian, the noncooperative elliptic system:
\be
-\Delta u=|u|^{p-1}u+f_u(x,u,v)
\ee
\be
\Delta v=|v|^{q-1}v+f_v(x,u,v)
\ee
\be
v\mid_{\partial Q}=u\mid_{\partial Q}=0.
\ee
However the approaches in \cite{CDHL},\cite{CastroClapp} are incomplete as they rely on estimating $\int_Q|\nabla{[\tau(u)u]}|^2dx$ for $u\in H_0^1(Q)$ but the functional $\tau:H_0^1(Q)\rightarrow \R$ is not Fr\'echet differentiable and the authors do not define what they mean by  $\nabla{[\tau(u)u]}$,$\int_Q|\nabla{[\tau(u)u]}|^2dx$, for arbitrary $u\in H_0^1(Q)$.}\\

In Section 1:There is a functional $I_\beta$ whose critical points correspond formally to solutions of (\ref{mpapier}). However as indicated by the approach of \cite{R82}, for technical reasons we will work with another functional $J_{\beta}$. We prove Palais-Smale conditions at large energies independently of $\beta$ for the functional $J_{\beta}$ and show implications for the functional $I_\beta$.


In Section 2: We construct the map $H$ whose energy is bounded independenlty of $\beta$. This is the main novelty of the paper which leads to the compactness needed to show the existence of classical solutions.

In Section 3 we adapt the arguments of \cite{R78} and \cite{R84} to end the proof. First we show that $u(\beta)$ is a classical solution of the modified equation (\ref{mpapier}) then we obtain a $C^0$ estimate for $w(\beta)$. This is followed by a $C^0$ on $v(\beta)$, and the existence of a $C^0$-solution $u$ is proved. We then use the bootstrapping argument in \cite{R78} to prove the existence of classical solutions. The multiplicity is deduced by noticing the lower estimates on the critical values $c_n^m(\delta)$ go to infinity as $n\rightarrow\infty$.\\
Functional $I_\beta$:\\
We define the functional $I_\beta $:
\be
I_\beta(u)=\int_Q[\frac{1}{2}(u_t^2-u_x^2-\beta v_t^2)-\frac{1}{s+1}|u|^{s+1}-fu]dxdt,
\ee
the domain $Q=[0,\pi]\times[0,2\pi]$.
We seek time-periodic solutions satisfying Dirichlet boundary conditions so we seek functions $u\in\R$ with expansions of the form
\[
u(x,t)=\sum_{(j,k)\in\N\times\Z}\widehat{u}(j,k)\sin jx e^{ikt}
\]
and define the function space
\[
||u||_{E^p}=\sum_{j\neq |k|}|Q||k^2-j^2|^s|\widehat{u}(j,k)|^2+\sum_{j=\pm k}|\widehat{u}(j,k)|^2
\]
where we denote by $E$ the space $E^p$ with $p=1$. Define the functions spaces $E^+,E^-,N$ as follows:
\[
N=\{ u\in E, \,\, \widehat{u}(j,k)=0 \,\, {\rm{for}} \,\, j\neq |k|\}
\]
\[
E^+=\{ u\in E, \,\, \widehat{u}(j,k)=0 \,\, {\rm{for}} \,\, |k|\leq j\}
\]
\[
E^-=\{ u\in E, \,\, \widehat{u}(j,k)=0 \,\, {\rm{for}} \,\, |k|\geq j\},
\]
$w=w^++w^-$ where $w^+\in E^+$,$w^-\in E^-$ and $v\in N$
and define the norm on $E\oplus N$
\[
||u||_{\beta,E}^2=||w^+||^2_E+||w^-||^2_E+\beta||v_t||^2_{L^2}.
\]
where 
\be
\rm{If }  \,\, j<|k|, \,\ w^+(j,k)=\sqrt{|Q|}\hat{u}(j,k).
\ee
\be
\rm{If }  \,\, j>|k|, \,\ w^-(j,k)=\sqrt{|Q|}\hat{u}(j,k).
\ee
\be
\rm{If }  \,\, j=|k|, \,\ v(j,k)=\sqrt{|Q|}\hat{u}(j,k).
\ee
When $u$ is trigonometric polynomial, $I_\beta$ can also be represented as:
\be
I_\beta(u)=\frac{1}{2}(||w^+||_E^2-||w^-||_E^2-\beta||v_t||_{L^2}^2)-\frac{1}{s+1}||u||_{L^{s+1}}^{s+1}-\int_Q fudxdt .
\ee
This is true because $\{\frac{\sin jx e^{ikt}}{\sqrt{2|Q|}},\frac{\cos jx e^{ikt}}{\sqrt{2|Q|}} \}$ is an orthonormal basis of $L^2(x,t)[0,2\pi]\times[0,2\pi]$ with periodic boundary conditions:
\[
u(x,t)=\sum_{(j,k)\in\N\times\Z}\widehat{u}(j,k)\sqrt{2|Q|}\frac{\sin jx e^{ikt}}{\sqrt{2|Q|}}
\]

Now $u$ being odd in $x$ can be extended to a doubly periodic functions and we will use Parseval formula to estimate the integrals of $u_t,u_x$ in Fourier basis.
\begin{eqnarray} \int_0^{\pi}\int_0^{2\pi}u_x^2dxdt
& = & \frac{1}{2}\int_0^{2\pi}\int_0^{2\pi}u_x^2dxdt\nonumber\\
& = & \frac{1}{2}\int_0^{2\pi}\int_0^{2\pi}|u_x|^2dxdt\nonumber\\
& = & \frac{1}{2}\sum_{j,k}j^22|Q||\hat{u}(j,k)|^2\nonumber\\
& = &\sum_{j,k}j^2|Q||\hat{u}(j,k)|^2
\end{eqnarray}
The first equality hold because $u$ is odd in $x$ and $u_x$ is even in $x$, and the second because $u$ is real. Then the second follows from Parseval's formula in $L^2(x,t)[0,2\pi]\times[0,2\pi]$.
Similarly 
\begin{eqnarray}
\int_0^{\pi}\int_0^{2\pi}u_t^2dxdt
& =& \int_0^{\pi}\int_0^{2\pi}|u_t|^2dxdt\nonumber\\
& = & \sum_{j,k}k^2|Q||\hat{u}(j,k)|^2
\end{eqnarray}
Since $u$ is real valued we can replace $u_t^2$ by its modulus $|u_t|^2$. Same earlier with $u_x$.
We can now conclude that
\begin{eqnarray}
\int_0^{\pi}\int_0^{2\pi}u_t^2-u_x^2dxdt
& = & \sum_{j,k}(k^2-j^2)|Q||\hat{u}(j,k)|^2\nonumber\\
& = & \sum_{k>j}|k^2-j^2|Q||\hat{u}(j,k)|^2-\sum_{k<j}|k^2-j^2||Q||\hat{u}(j,k)|^2\nonumber\\
& = & ||w^+||_E^2-||w^-||_E^2
\end{eqnarray}
The spectrum of the linear operator $\partial_t^2-\partial_x^2$ under Dirichlet boundary conditions in space and time-periodicity consists of
\[
-k^2+j^2
\]
where the eigenfunctions are the $\sin jx \cos kt,\sin jx\sin kt$. The eigenfunctions here are ordered as in \cite{Tanaka88} i.e
\[
...-\mu_3\leq-\mu_2\leq-\mu_1<0<\mu_1\leq\mu_2\leq\mu_3\leq...
\]
where the $\mu_l$ are the eigenvalues of $\partial_t^2-\partial_x^2$ and have multiplicity one. Rearranging the eigenvalues this way is possible because all the non-zero eigenspaces of $\partial_t^2-\partial_x^2$ have finite multiplicity. The $\mu_l\rightarrow+\infty$ as $l\rightarrow+\infty$  and denote by $e_l$ the corresponding eigenfunctions,
and we define the spaces
\[
{E^{+n}}=span\{e_l, 1\leq l\leq n \}.
\]
 For the Galerkin procedure we define the spaces
\[
E^{-m}=span\{\sin jx\cos kt, \sin jx\sin kt, \,\,\, j+k\leq m \,\ j<k \},
\]
\[
N^m=span\{ \sin jx\cos jt, \sin jx\sin jt, \,\,\ j\leq m \}
\]
which are employed in the minimax procedure.

We start by following the procedure of \cite{R82} for perturbation problems by proving some properties of the functional $I_\beta$. The difference here is that additionally we show that the constants involved in all the proof are independent of $\beta$ to prepare for passing to the limit as $\beta\rightarrow 0$.
\begin{lem} Suppose that $u$ is a critical point of $I_\beta$. Then there is a constant $a_6$ depending on $s,f$ but independent of $\beta$ such that
\be
\int_Q\frac{|u|^{s+1}}{s+1}dxdt\leq a_6[I_\beta^2(u)+1]^{\frac{1}{2}}
\label{us+1}
\ee
\label{1.1}\end{lem}
Proof:
\begin{eqnarray}I_\beta(u)
& = & I_\beta(u)-\frac{1}{2}I_\beta^\prime(u)u \nonumber\\
& = & \frac{s-1}{2(s+1)}\int_Q|u|^{s+1}dxdt-\frac{1}{2}\int_Qfudxdt.\nonumber\\
& \geq & \frac{s-1}{4(s+1)}\int_Q|u|^{s+1}dxdt-\int_Qfudxdt.
\label{fu}\end{eqnarray}
Now recalling Young inequality
\be
ab\leq \frac{a^p}{p}+\frac{b^q}{q}
\ee
for $\frac{1}{p}+\frac{1}{q}=1,a,b>0$ to $\int_Q\frac{1}{\epsilon}f\epsilon udxdt$ we deduce for $\epsilon>0$ small we have

\be
I_\beta(u)\geq \frac{s-1}{4(s+1)}\int_Q|u|^{s+1}dxdt-\frac{1}{\epsilon}\frac{s}{s+1}||f||_{L^{\frac{s+1}{s}}}^{\frac{s+1}{s}}-\epsilon\frac{1}{s+1}||u||_{L^{s+1}}^{s+1}
\ee
\be
\frac{1}{\epsilon}\frac{s}{s+1}||f||_{L^{\frac{s+1}{s}}}^{\frac{s+1}{s}}+I_\beta(u)\geq (\frac{s-1}{4(s+1)}-\epsilon\frac{1}{s+1})\int_Q|u|^{s+1}dxdt
\ee
\be
\max(1,\frac{1}{\epsilon}\frac{s}{s+1}||f||_{L^{\frac{s+1}{s}}}^{\frac{s+1}{s}})(1+I_\beta(u))\geq (\frac{s-1}{8(s+1)})\int_Q|u|^{s+1}dxdt
\ee
for $\epsilon(s)<< 1$ small enough and independent of $\beta$ hence
\be
\int_Q|u|^{s+1}dxdt\leq\frac{8(s+1)}{s-1} \frac{1}{\epsilon}\frac{s}{s+1}||f||_{L^{\frac{s+1}{s}}}^{\frac{s+1}{s}}(1+I_\beta(u))
\ee
and
\be
\int_Q\frac{|u|^{s+1}}{s+1}dxdt\leq\frac{4\sqrt{2}}{s-1} \frac{1}{\epsilon}\frac{s}{s+1}||f||_{L^{\frac{s+1}{s}}}^{\frac{s+1}{s}}(1+I^2_\beta(u))^{\frac{1}{2}}
\ee
 while we choose $a_6=\frac{4\sqrt{2}}{s-1} \frac{1}{\epsilon}\frac{s}{s+1}||f||_{L^{\frac{s+1}{s}}}^{\frac{s+1}{s}}$.\\
We define the functional $J_{\beta}$ which is amenable to minimax procedure. We start by defining a bump function $\chi$. $\chi\in C^\infty(\R,\R)$:
\begin{equation}
\left\{ \begin{array}{ll}   \chi(t)=1 ,\,\ \rm{if} \,\ t\leq 1  & \\
\chi(t)=0 \,\, \rm{if} \,\, t> 2 \,\, . &
\end{array} \right.
\end{equation}
and $-2<\chi^\prime<0$, for $1<t<2$. Then define
\be
{\cal{I}}_\beta(u)=2a_6(I_\beta^2(u)+1)^{\frac{1}{2}}
\ee
and
\[
\psi(u)=\chi({\cal{I}}_\beta^{-1}(u)\int_Q\frac{|u|^{s+1}}{s+1}dxdt)
\]
\be
J_{\beta}(u)=\int_Q[\frac{1}{2}(u_t^2-u_x^2-\beta v_t^2)-\frac{1}{s+1}|u|^{s+1}-\psi(u)fu]dxdt,
\ee
which on $E^{+m}\oplus E^{-m}\oplus N^m$ can be rewritten as
\be
J_{\beta}(u)=\frac{1}{2}(||w^+||_E^2-||w^-||_E^2-\beta||v_t||_{L^2}^2)-\frac{1}{s+1}||u||_{L^{s+1}}^{s+1}-\int_Q \psi(u)fudxdt .
\ee
 
\begin{lem}{\rm{If}} 
\be
{\cal{I}}^{-1}_\beta(u)\int_Q\frac{1}{s+1}|u|^{s+1}dxdt\leq 2
\label{supp}\ee
{\rm{then is a constant}} $\alpha_3$ {\rm{independent of}} $\beta$ {\rm{such that}}
\[
|\int_Qfudxdt|\leq \alpha_3({I}_\beta^{\frac{1}{s+1}}(u)+1)
\]
\label{lemma1.13}
\end{lem}
Proof:
\[
|\int_Qfudxdt|\leq c(f,s)||u||_{L^{s+1}}
\]
by Holder inequality, then  since  we assumed (\ref{supp})
we have
\[
|\int_Qfudxdt|\leq c(f,s)||u||_{L^{s+1}}\leq\alpha_3(I_\beta^{\frac{1}{s+1}}(u)+1)
\]
where $\alpha_3$ depends on $f,s$
\begin{lem}{\rm{There is a constant}}  $\gamma_1$ {\rm{depending on}} $f,s$ {\rm{but independent of}} $\beta$ {\rm{such that}}
\be
|J_{\beta}(u)-J_{\beta}(-u)|\leq \gamma_1(|J_{\beta}(u)|^{\frac{s}{s+1}}+1)
\label{even}\ee
\end{lem}
Proof:\\
\be
J_{\beta}(u)-J_{\beta}(-u)=-\psi(u)\int_Qfudxdt+\psi(-u)\int_Qfudxdt
\label{Jb}\ee
Step1: If $\psi(u)\neq 0$ then
\be
|\int fudxdt|\leq \alpha_3(I_\beta^{\frac{1}{s+1}}(u)+1).
\label{phi}\ee
We have
\be
I_{\beta}(u)=J_{\beta}(u)-\int_Qfudxdt+\int_Q\psi(u)fudxdt
\ee
thus
\be
|I_\beta(u)|\leq |J_{\beta} (u)|+2|\int_Qfudxdt|
\label{IJ1}\ee
\be
|I_\beta(u)|^{\frac{1}{s+1}} \leq |J_\beta(u)|^{\frac{1}{s+1}}+2^{\frac{1}{s+1}}|\int fudxdt|^{\frac{1}{s+1}}
\label{IJ2}
\ee
\be
|\int_Qfudxdt| \leq \alpha_3(|J_\beta(u)|^{\frac{1}{s+1}}+2^{\frac{1}{s+1}}|\int_Qfudxdt|^{\frac{1}{s+1}}+1) 
\ee
Now
\be
\frac{x}{2}>\frac{A}{2}x^{\frac{1}{s+1}} \,{\rm{
 if}} \,\,  x>A^{\frac{s+1}{s}}
\ee
so if we choose $\frac{A}{2}=\alpha_32^{\frac{1}{s+1}}$ and $x=|\int_Qfudxdt|$, for 
\be
|\int_Qfudxdt|>(2\alpha_32^{\frac{1}{s+1}})^{\frac{s+1}{s}}\ee
we have
\be
\frac{1}{2}|\int_Qfudxdt|\leq |\int_Qfudxdt|-\alpha_32^{\frac{1}{s+1}}|\int_Qfudxdt|^{\frac{1}{s+1}}\leq\alpha_3(|J_\beta(u)^{\frac{1}{s+1}}+1)
\ee
hence
\be
|\int_Qfudxdt|\leq2\alpha_3(|J_\beta(u)^{\frac{1}{s+1}}+1)
\ee
and recalling (\ref{Jb}) we deduce
\[
|J_{\beta}(u)-J_{\beta}(-u)|\leq 4\alpha_3(|J_\beta(u)^{\frac{1}{s+1}}+1).
\]
If on the other hand 
\[
|\int_Qfudxdt|<(2\alpha_32^{\frac{1}{s+1}})^{\frac{1}{s+1}}
\]
then (\ref{Jb}) and $0\leq \phi(u),\psi(-u)\leq 1$ implies (\ref{even}).\\
Now Step 2: if $\phi(-u)\neq 0$\\ then the argument in the previous lemma \ref{lemma1.13} implies
\be
|\int f(-u)dxdt|\leq \alpha_3(I_\beta^{\frac{1}{s+1}}(-u)+1).
\ee
We also have
\be
I(u)-I(-u)=-2\int fudxdt
\ee
thus
\be
|\int f(-u)dxdt|\leq \alpha_3(I_\beta^{\frac{1}{s+1}}(u)+2^{\frac{1}{s+1}}|\int f(-u)dxdt|^{\frac{1}{s+1}}+1)
\ee
and by repeating the argument made earlier with $x=|\int_Qfudxdt|$ and $\frac{A}{2}=\alpha_32^{\frac{1}{s+1}}$, if we have 
\be
|\int_Qfudxdt|>(2\alpha_32^{\frac{1}{s+1}})^{\frac{s+1}{s}}
\label{A}\ee
we deduce
\be
|\int_Qfudxdt|\leq2\alpha_3(|I_\beta(u)^{\frac{1}{s+1}}+1).
\ee
Now 

Now recalling (\ref{IJ1},\ref{IJ2})
and (\ref{phi}) implies
\[
|\int_Qfudxdt|\leq 2\alpha_3(|J_\beta(u)^{\frac{1}{s+1}}+2^{\frac{1}{s+1}}|\int_Qfudxdt|^{\frac{1}{s+1}}+1)
\]
\be
|\int_Qfudxdt|-2\alpha_32^{\frac{1}{s+1}}|\int_Qfudxdt|^{\frac{1}{s+1}}\leq 2\alpha_3(|J_\beta(u)^{\frac{1}{s+1}}+1)
\ee
now the inequality (\ref{A}) implies with an argument similar to the given earlier in Step 1, that
\be
\frac{1}{2}|\int fudxdt|>\alpha_32^{\frac{1}{s+1}}|\int fudxdt|^{\frac{1}{s+1}}
\ee
and we conclude again.
\be
|\int_Qfudxdt|\leq2(2\alpha_3)(|J_\beta(u)^{\frac{1}{s+1}}+1).
\ee
If $|\int_Qfudxdt|<(2\alpha_32^{\frac{1}{s+1}})^{\frac{s+1}{s}}$ the lemma follows again.

Step 3:\\

If $\psi(u)=\psi(-u)=0$ then by (\ref{Jb}) the lemma follows again.


\begin{lem}{\rm{There are constants}} $\alpha_0,M_0>0$ {\rm{depending on}} $f,s$ {\rm{independent of}} $\beta$ {\rm{such that whenever}} $M\geq M_0$, {\rm{then}} $J_\beta(u)\geq M$
{\rm{and}}  when (\ref{supp}) {\rm{is satisfied}} 
  {\rm{then}} $I_\beta(u)\geq \alpha M_0$
\end{lem}
Proof:\\
\be
I_\beta(u)\geq J_\beta(u)-2|\int_Qfudxdt|
\label{ineq1}\ee
while if (\ref{supp}) {\rm{is satisfied}} then
 there is $\alpha_3$ independent of $\beta$ such that 
\be
2\alpha_3(|I_\beta(u)|^{\frac{1}{s+1}}+1)\geq 2|\int_Qfudxdt|
\label{ineq0}\ee
or
\be
|I_\beta(u)|^{\frac{1}{s+1}}\geq \frac{1}{\alpha_3}|\int_Qfudxdt|-1
\label{ineq2}\ee
and adding (\ref{ineq1}) and (\ref{ineq0})
\be
I_\beta(u)+2\alpha_3|I_\beta(u)|^{\frac{1}{s+1}}\geq J_\beta(u)-C\geq \frac{M}{2}
\label{ineq3}
\ee
for $M_0$ large enough. If $I_\beta(u)\leq 0$, then by Young inequality
\be
2\alpha_3|I_\beta(u)|^{\frac{1}{s+1}}\leq \frac{(2\alpha_3)^{\frac{s+1}{s}}}{\frac{s+1}{s}}+\frac{1}{s+1}|I_\beta(u)|
\ee
while the inequality (\ref{ineq3}) implies
\be
2\alpha_3|I_\beta(u)|^{\frac{1}{s+1}}\geq -I_\beta(u)+\frac{M}{2}
\ee
hence
\be
\frac{(2\alpha_3)^{\frac{s+1}{s}}}{\frac{s+1}{s}}+\frac{1}{s+1}|I_\beta(u)|\geq-I_\beta(u)+\frac{M}{2} =|I_\beta(u)|+\frac{M}{2}
\ee
thus 
\be
\frac{s}{s+1}|I_\beta(u)|\leq -\frac{M}{4}<0
\ee
and we have a contradiction.
\begin{lem}Lemma 1.29 \cite{R82} In $E^{+m}\oplus E^{-m}\oplus N^m$,there is a constant $M_1>0$ independent of $\beta,m$ such that $J_\beta(u)\geq M_1$ and $J_\beta^\prime(u)=0$ implies that
$J_\beta(u)=I_\beta(u)$ and $I_\beta^\prime(u)=0$
\end{lem}
Proof:\\
We follow step by step the argument in \cite{R82}.\\ It suffices to show that
\be
{\cal{I}_{\beta}}^{-1}(u)\int_Q\frac{1}{s+1}|u|^{s+1}dxdt\leq 1
\ee
\be
J_\beta^\prime(u)u=\int_Qw_t^2-w_x^2-\beta v_t^2-|u|^{s+1}dxdt-\psi(u)\int_Qfudxdt-\psi^\prime(u)u\int_Qfudxdt
\ee
where
\begin{eqnarray} \psi^\prime(u)u
& = & \chi^\prime({\cal{I}_{\beta}}^{-1}(u)\int_Q\frac{1}{s+1}|u|^{s+1}dxdt) \nonumber\\
&   & \times [-{\cal{I}_\beta}^{-3}(u)(2a_6)^22I_\beta(u)I^\prime_\beta(u)u\int_Q\frac{|u|^{s+1}}{s+1}dxdt
      +{{\cal{I}_\beta}}^{-1}(u)\int_Q|u|^{s+1}dxdt]\nonumber
\end{eqnarray}
and
\be
J^\prime_\beta(u)u=(1+T_1(u))\int_Qw_t^2-w_x^2-\beta v_t^2dxdt-(1+T_2(u))\int_Q|u|^{s+1}dxdt-(\psi(u)+T_1(u))\int_Qfudxdt
\ee
where $T_1,T_2$ are exactly as in \cite{R82}:
\be
T_1(u)=\chi^\prime({\cal{I}_{\beta}}^{-1}(u)\int_Q\frac{1}{s+1}|u|^{s+1}dxdt)(2a_6)^2{\cal{I}_\beta}^{-3}(u)\int_Q\frac{|u|^{s+1}}{s+1}dxdt\int_Qfudxdt
\ee
and
\be
T_2(u)=\chi^\prime({\cal{I}_{\beta}}^{-1}(u)\int_Q\frac{1}{s+1}|u|^{s+1}dxdt){\cal{I}_\beta}^{-1}(u)\int_Qfudxdt+T_1(u)
\ee
and the conclusion follows just as in \cite{R82}.

\begin{eqnarray}J_\beta(u)-\frac{1}{2(1+T_1)}J^\prime_\beta(u)u
& = & \int_Q(-\frac{1}{s+1}+\frac{1+T_2}{2(1+T_1)})|u|^{s+1}dxdt-[\psi(u)-\frac{\psi(u)+T_1}{2(1+T_1)}]\int_Qfudxdt\nonumber
\end{eqnarray}
Now as $T_1,T_2\rightarrow 0$ as $M\rightarrow +\infty$, if $u$ is a critical point of $J_\beta$ we have:
\begin{eqnarray}I_\beta(u)-\psi(u)\int fudxdt+\int fudxdt
& = & \int_Q(-\frac{1}{s+1}+\frac{1+T_2}{2(1+T_1)})|u|^{s+1}dxdt\nonumber\\
&   &  -[\psi(u)-\frac{\psi(u)+T_1}{2(1+T_1)}]\int_Qfudxdt\nonumber
\end{eqnarray}
\begin{eqnarray}I_\beta(u)+\int fudxdt
& = & \int_Q(-\frac{1}{s+1}+\frac{1+T_2}{2(1+T_1)})|u|^{s+1}dxdt\nonumber\\
&   &  +[\frac{\psi(u)+T_1}{2(1+T_1)}]\int_Qfudxdt\nonumber
\end{eqnarray}
\begin{eqnarray}I_\beta(u)
& = & \int_Q(-\frac{1}{s+1}+\frac{1+T_2}{2(1+T_1)})|u|^{s+1}dxdt \nonumber\\
&   &  +[\frac{\psi(u)}{2(1+T_1)}+\frac{T_1}{2(1+T_1)}-1]\int_Qfu dxdt\nonumber
\end{eqnarray}
Now as $\frac{1+T_2}{2(1+T_1)}\rightarrow \frac{1}{2}$,$\frac{T_1}{2(1+T_1)}\rightarrow 0$,$-\frac{1}{s+1}+\frac{1+T_2}{2(1+T_1)}\rightarrow \frac{s-1}{2(s+1)},0\leq \psi(u)\leq 1$ and 
\begin{eqnarray}I_\beta(u)
& \geq & \int_Q\frac{s-1}{4(s+1)}|u|^{s+1}dxdt -\frac{7}{4}\int_Qfu dxdt\nonumber
\end{eqnarray}
The preceding inequality is essentially (\ref{fu}) where $\frac{7}{4}$ replaces $1$. Now following step by step the computations in lemma \ref{1.1} we obtain
\be
\int_Q\frac{|u|^{s+1}}{s+1}dxdt\leq 2a_6(1+I_\beta^2)^{\frac{1}{2}}
\ee

We now show that the functional $J_\beta$ satisfies the Palais-Smale condition at large energies independently of $\beta$ in  $E^{+m}\oplus E^{-m}\oplus N^m$.
\begin{lem}{\rm{There is a constant}} $M_2$ {\rm{independent of}} $\beta$ {\rm{such that the Palais-Smale condition is satisfied on}} $A_{M_2}=\{ u\in E^{+m}\oplus E^{-m}\oplus N^m, \,\,\ J_\beta(u)\geq M_2 \}$.
\end{lem}
Proof:\\
Let $u_l=w_l+v_l=w_l^++w_l^-+v_l$ a Palais-Smale sequence at large energies, there are $M_2,K$ independent of $\beta,m$ such that if $M_2\leq J_\beta(u_l)\leq K$
and $J_\beta^\prime(u_l)\rightarrow 0$, then $||u_l||_{E,\beta}\leq c(\beta)$ hence, since $E^{+m}\oplus E^{-m}\oplus N^m$ is finite dimensional, $u_l$ has a convergent subsequence.\\ 
\begin{eqnarray}J_\beta(u_l)-\rho J_\beta^\prime(u_l)(u_l)
& = & (\frac{1}{2}-\rho(1+T_1(u_l))) \int_Qw_{lt}^2-w_{lx}^2-\beta v_{lt}^2dxdt\nonumber\\
&   & +[\rho(1+T_2(u_l))-\frac{1}{s+1}]\int_Q|u_l|^{s+1}dxdt\nonumber\\
&   &    +(\rho(\psi(u_l)+T_1(u_l))-\psi(u_l))\int_Qfu_ldxdt
\end{eqnarray}
now we choose $\rho=\frac{1}{2(1+T_1(u_l))}$ then we have
\[
\rho\rightarrow\frac{1}{2} \,\,{\rm{independently}} \,\ {\rm{of}} \,\, \beta \,\ {\rm{as}} \,\ M_2\rightarrow+\infty
\]
\begin{eqnarray}J_\beta(u_l)-\rho J_\beta^\prime(u_l)(u_l)
& = & [\rho(1+T_2(u_l))-\frac{1}{s+1})]\int_Q|u_l|^{s+1}dxdt\nonumber\\
       +(\rho(\psi(u_l)+T_1(u_l))-\psi(u_l))\int_Qfu_ldxdt\nonumber\\
& \geq & [\rho(1+T_2(u_l))-\frac{1}{s+1})-\frac{\epsilon(s)}{s+1}]\int_Q|u_l|^{s+1}dxdt-c(f,s)\nonumber
\end{eqnarray}
where $\epsilon(s)$ can be chosen to be a small positive constant by applying Young inequality so that
\be
 [\rho(1+T_2(u_l)-\frac{1}{s+1})-\frac{\epsilon(s)}{s+1}]>0,
\ee and $c(f,s)$ is another constant depending on $f,s$, both being independent of $\beta$. Now recall that $J^\prime_\beta(u_l)\rightarrow 0$ and $\rho\rightarrow \frac{1}{2}$
\be
J_\beta(u_l)-\rho J^\prime_\beta(u_l)u_l\leq K+\rho||u_l||_{E,\beta}
\ee
so we have the inequalities:
\be
K+\rho||u_l||_{\beta,E}\geq J_\beta(u_l)-\rho J^\prime(u_l)u_l\geq c_3(s)||u||_{L^{s+1}}^{s+1}-c_2(f,s)
\ee
thus
\be
\int_Q|u_l|^{s+1}dxdt\leq c_4(f,s)||u_l||_{E,\beta}+K+c_2(f,s).
\label{ineq5}\ee
Now
\be
J_\beta^\prime(u_l)v_l=(1+T_1(u_l))\int_Q\beta v_{lt}^2dxdt-(1+T_2(u_l))\int_Q|u_l|^{s-1}u_lv_ldxdt-(\psi(u_l)+T_1(u_l))\int_Qfv_ldxdt.
\ee
$u_l$ is a Palais-Smale sequence so there exists $\epsilon$ small such that \\
$J_\beta^\prime(u_l)v_l\leq\epsilon||v_l||_{\beta,E}$ thus
\[
(1+T_1(u_l))\beta||v_{lt}||_{L^2}^2\leq (1+T_2(u_l))\int_Q|u_l|^{s-1}u_lv_ldxdt+(\psi(u_l)+T_1(u_l))\int_Qfv_ldxdt+\epsilon||v_l||_{\beta,E}.
\]
Now for $M_2$ large enough (independently of $\beta$) and we have
\be
\frac{1}{2}\beta||v_{lt}||_{L^2}^2\leq (2\int_Q|u_l|^{s}|v_l|dxdt+2\int_Q|f||v_l|dxdt+\epsilon||v_l||_{\beta,E}
\ee
 and applying H\"older inequality we deduce:
\[
\beta||v_{lt}||_{L^2}^2\leq 4||u_l||_{L^{s+1}}^s||v_l||_{L^{s+1}}+4||v_l||_{L^{s+1}}||f||_{L^{\frac{s+1}{s}}}+2\epsilon ||v_l||_{\beta,E}.
\]
 A similar computation gives
\be
||w_l^+||_{E,\beta}^2\leq 4||u_l||_{L^{s+1}}^s||w_l^+||_{L^{s+1}}+4||w_l^+||_{L^{s+1}}||f||_{L^{\frac{s+1}{s}}}+2\epsilon||w_l^+||_{E},
\ee
and
\be
||w_l^-||_{E,\beta}^2\leq 4||u_l||_{L^{s+1}}^s||w_l^-||_{L^{s+1}}+4||w_l^-||_{L^{s+1}}||f||_{L^{\frac{s+1}{s}}}+2\epsilon||w_l^-||_{E}.
\ee
\begin{eqnarray}||u_l||_{E,\beta}^2
& \leq & (4||u||^s_{L^{s+1}}+4||f||_{L^{\frac{s+1}{s}}})(||v_l||_{L^{s+1}}+||w_l^+||_{L^{s+1}}+||w_l^-||_{L^{s+1}})\nonumber\\
&   & +2\epsilon(||v_l||_{E,\beta}+||w_l^+||_{E}+||w_l^-||_{E})
\end{eqnarray}

We now estimate $||v_l||_{L^{s+1}}$:$v_l=u_l-w_l^+-w_l^-$ hence
\begin{eqnarray}||v_l||_{L^{s+1}}
& \leq & ||u_l||_{L^{s+1}}+||w_l^+||_{L^{s+1}}+||w_l^-||_{L^{s+1}}\nonumber\\
& \leq & c_5||u_l||_{E,\beta}^{\frac{1}{s+1}}+c_6+c_6||w_l^-||_E+c_6||w_l^+||_E
\label{ineq4}\\
& \leq & c_7||u_l||_{E,\beta}+c_8
\end{eqnarray}
where the constants $c_5,c_6,c_7,c_8$ depend on $s,K,c_4,c_2D(f,s)$ are independent of $\beta$ and (\ref{ineq4}) follows from (\ref{ineq5}) and the Sobolev inequality $||w_l||_{L^p}\leq c(p)||w_l||_E$,
Now from 
\begin{eqnarray}||u||_{L^{s+1}}
& \leq & c_5||u||_{E,\beta}^{\frac{1}{s+1}}+c_6
\end{eqnarray}
we deduce
\begin{eqnarray}||u||^s_{L^{s+1}}
& \leq & 2^sc_5^s||u||_{E,\beta}^{\frac{s}{s+1}}+2^sc_6^s
\end{eqnarray}
 We can now deduce:
\begin{eqnarray}||u_l||_{E,\beta}^2
& \leq & c(f)(1+2^sc_5^s||u||_{E,\beta}^{\frac{s}{s+1}}+2^sc_6^s)(c_7||u_l||_{E,\beta}+c_8+2c_6||u_l||_{E,\beta})+c||u_l||_{E,\beta}\nonumber\\
& \leq & c||u||_{E,\beta}^{\frac{2s+1}{s+1}}+d
\end{eqnarray}
where $c,d$ depend on the constants above but independent of $\beta$. Since \\
${\frac{2s+1}{s+1}}<2$, we deduce $||u_l||_{E,\beta}<+\infty$. $E^{+m}\oplus E^{-m}\oplus N^m$ is finite dimensional, so Palais-Smale is satisfied.
\section{Estimates on minimax values independently of $\beta$}
\begin{lem} There is $R_n\rightarrow +\infty$ such that $J_\beta(u)\rightarrow -\infty$, uniformly as $||u||_{\beta,E}=R_n\rightarrow +\infty$, {\rm{for}} $u\in E^{+n}\oplus E^{-m}\oplus N^m$. As a result we can also assume that $R_{n+1}>4R_n$, without loss of generality.
\label{R_n}
\end{lem}
Proof:\\
Let $u=w^++w^-+v$
\be
||u||_{E,\beta}^2=R_n^2
\ee then either $||w^+||^2_{E,\beta}\geq \frac{R_n^2}{3}$ or $||w^-|||^2_{E,\beta}+\beta||v_t||_{L^2}^2\geq \frac{2R_n^2}{3}$.\\
Case 1: $||w^+||_E^2\geq \frac{R_n^2}{3}$:\\
\begin{eqnarray}J_\beta(u)
& = & \frac{1}{2}(||w^+||_E^2-||w^-||_E^2-\beta||v_t||^2_{L^2})-\frac{1}{s+1}||u||_{L^{s+1}}^{s+1}- \psi(u)\int_Q fudxdt\nonumber\\
& \leq & \frac{1}{2}||w^+||_E^2-\frac{1}{s+1}||u||_{L^{s+1}}^{s+1}- \psi(u)\int_Q fudxdt\nonumber\\
& \leq & \frac{1}{2}||w^+||_E^2-a(s)||u||_{L^{s+1}}^{s+1}+c(f)
\end{eqnarray}
where $a(s),c(f)>0$ by Young inequality and as $w^+\in E^{+n}$ and $s>1$ we also have:
\be
\frac{1}{\mu_n}||w^+||_E\leq||w^+||_{L^2}\leq||u||_{L^2}\leq c(Q)
||u||_{L^{s+1}}
\ee
thus
\begin{eqnarray}J_\beta(u)
& \leq &  \frac{1}{2}||w^+||_E^2-a(s)(\frac{||w^+||_E}{c(Q)\mu_n})^{s+1}+c(f)
\end{eqnarray}
and for $R_n(s,f,\mu_n)$ large enough,$J_\beta(u)\rightarrow -\infty$ uniformly as $s>1$.\\
Case 2:$||w^+||^2< \frac{R_n^2}{3}$ hence $||w^-||_E^2+\beta||v_t||_{L^2}^2\geq \frac{2R_n^2}{3}$, thus:
\begin{eqnarray}J_\beta(u)
& \leq &-\frac{R_n^2}{6}-\frac{1}{s+1}||u||_{L^{s+1}}^{s+1}-\psi(u)\int_Q fudxdt\nonumber\\
& \leq & -\frac{R_n^2}{6}-a(s)||u||_{L^{s+1}}^{s+1}+c(f),
\end{eqnarray}
by Young inequality, and we can conclude again that $R_n$ large enough \\
$J_\beta(u)\rightarrow-\infty$ uniformly which ends the proof of the lemma.
We now define the minimax sets and the minimax values which will lead to the existence of critical values:\\
Let $B(R,W)$ the closed ball, of radius $R$, in  a subspace $W$ of $E\oplus N$:
\[
B(R,W)=\{u\in W,\,\, ||u||_{E,\beta}\leq R\}
\]
\[
D_n^m=\{u\in E^{+n}\oplus E^{-m}\oplus N^{m} \,\, {\rm{and}} \,\, ||u||_{E,\beta} \leq R_n \}
\]
\[
\Gamma_n^m=\{ h:D_n^m \rightarrow E^{+m}\oplus E^{-m}\oplus N^m,\,\ h \, {\rm{continuous \,and \,odd}} \,\ , h(u)=u \, {\rm{for}}\, ||u||_{E,\beta}=R_n \}
\]
\be
b_n^m=\inf_{h\in\Gamma_n^m}\max_{u\in D_n^m}J_\beta(h(u))
\label{bnm}
\ee
\[
U_n^m=\{ u=t e_{n+1}+w, \,\, t\in[0,R_{n+1}], w\in B(R_{n+1},E^{+n}\oplus E^{-m}\oplus N^{m}), ||u||_{E,\beta}\leq R_{n+1} \}
\]
\begin{displaymath}
\Lambda_n^m=\left\{ \begin{array}{l} H\in C(U_n^m,E^{+m}\oplus E^{-m}\oplus N^{m}),\,\, \\
H_{\mid D_n^m}\in\Gamma_n^m,\, {\rm{and}} \,\, H(u)=u \,{\rm{if}} \,||u||_{E,\beta}=R_{n+1}, \,{\rm{ or \, if }} \\
 u\in B(R_{n+1},E^{+n}\oplus E^{-m}\oplus N^{m})\setminus B({R_n},E^{+n}\oplus E^{-m}\oplus N^{m})\end{array}  \right\}
\end{displaymath}
where the constants $R_n$ do not depend on $\beta,m$.
\[
\Lambda_n^m(\delta)=\{ H\in \Lambda_n^m, \,\, J_\beta(H(u))\leq b_n^m+\delta \,\, {\rm{on}} \,\, D_n^m \, \}
\]
\[
c_n^m=\inf_{H\in\Lambda_n^m}\max_{u\in U_n^m}J_\beta(H(u))
\]
and
\[
c_n^m(\delta)=\inf_{H\in\Lambda_n^m(\delta)}\max_{u\in U_n^m}J_\beta(H(u))
\]



\begin{lem}$\forall u\in D_n^m\cap E^{+n}$, there is a constant $C(n)$ independent of $\beta,m$ such that
\be
J_\beta(u)\leq C(n).
\ee
\label{D_n^m}
\end{lem}
Proof:\\
Let $u\in E^{+n}$,$u=w^++w^-+v$ where $w^-=v=0$,
\begin{eqnarray}J_\beta(u)
& =    & \frac{1}{2}||w^+||_E^2-\frac{1}{2}||w^-||_E^2-\frac{1}{2}\beta||v_t||_{L^2}^2-\int_Q\frac{|u|^{s+1}}{s+1}dxdt-\psi(u)\int_Qfudxdt\nonumber\\
& = & \frac{1}{2}||w^+||_E^2-\int_Q\frac{|w^+|^{s+1}}{s+1}dxdt-\psi(w^+)\int_Qfw^+dxdt\nonumber\\
& \leq & \frac{1}{2}||w^+||_E^2-\frac{1}{2}\int_Q\frac{|u|^{s+1}}{s+1}dxdt+c(f,s)\\
& \leq & c(f,s)+\sup_{u\in E^{+n}}\frac{1}{2}||w^+||_E^2-\frac{1}{2}\int_Q\frac{|w^+|^{s+1}}{s+1}dxdt\nonumber\\
& \leq & c(f,s)+\sup_{w^+\in E^{+n}}\frac{1}{2}||w^+||_E^2-c(s,Q)||w^+||_{L^2}^{s+1}.
\end{eqnarray}
Now in $E^{+n}$
\be
||w^+||_E^2\leq\mu_n||w^+||_{L^2}^2
\ee
and on the other-hand
\be
\sup_{w^+\in E^{+n}}\frac{1}{2}||w^+||_E^2-c(s,Q)||w^+||_{L^2}^{s+1}>0
\ee
as $s>1$, and is attained at say $\overline{u}$,
hence we have
\be
c(s,Q)||\overline{u}||_{L^2}^{s+1}\leq \frac{1}{2}||\overline{u}||_E^2\leq \frac{1}{2}\mu_n||\overline{u}||_{L^2}^2
\ee
and we can conclude there is $C(n)$ depending on $n$ but independent of $\beta$ such that
\be
J_{\beta}(u)\leq C(n).
\label{JbetaGamma}\ee
We now construct the map which leads to upper estimates independently of $\beta$ which is the main contribution of the paper.
\begin{theo}Let $0<\delta<c_n^m-b_n^m$, then
and there is a map $H\in \Lambda_n^m(\delta)$ such that
\be
J_\beta(H(u,t))\leq C(n+1)
\ee
in $U_n^m$ where $C(n+1)$, is independent of $\beta,m$.
\end{theo}
Proof:\\
Let $h\in\Gamma_n^m$ a minimizing map for $b_n^m$, (\ref{bnm}), i.e.:
\be
J_\beta(h(u))\leq b_n^m+\frac{\delta}{2}
\label{minim}
\ee
on $D_n^m$.\\
The aim is to construct a function $H(u,t)$ which is the identity map when $||u||_{E,\beta}=R_{n+1}$ and which coincides with a map $h$ at $t=0$, for which $J_\beta(H(u,t))\leq c(n+1)$ a constant independent of the small parameter $\beta$.\\
Let $u\in E^{+n}$ we have $J_\beta(h(u))\leq b_m^n+\frac{\delta}{2}$ where $b_m^n$ is bounded independently of $m,\beta$.  The idea is to deform $H(u,t)$ from $h$ at $t=0$ to a map whose range is in $E^{n+1}$ and then to the identity map while keeping $H(u,t)=u+te_{n+1}$, when $||u||_{E,\beta}=R_{n+1}.$$h$ also satisfies $h(0)=0$ which plays an important role in the proof.\\

$0\leq t\leq 1 $: We construct a map $H(u,t)$ whose range is in $E^{+n+1}$ in $R_n<||u||_{E,\beta}< 3R_n$:\\
\[
\left \{\begin{array}{cc}
0\leq ||u||_{E,\beta}\leq R_n &  H(u,t)=H_{1,1}(u,t)=h(u) \\
R_n\leq ||u||_{E,\beta}\leq 2R_n &  H(u,t)=H_{1,2}(u,t)=(1-t)u+t(-u+\frac{2R_n}{||u||_{E,\beta}}u) \\
2R_n\leq ||u||_{E,\beta}\leq 3R_n & H(u,t)=H_{1,3}(u,t)=(1-t)u+t(3u-\frac{6R_nu}{||u||_{E,\beta}})+t(\frac{||u||_{E,\beta}}{R_n}-2)e_{n+1} \\
3R_n\leq ||u||_{E,\beta}\leq R_{n+1} & H(u,t)=H_{1,4}(u,t)=u+te_{n+1} \end{array} \right.
\]
 We now have to show the continuity of the function $H$ constructed. Since $h\in \Gamma_n^m$, $h(u)=u$ for $||u||_{E,\beta}\geq R_n$, we have to verify that gluing patches of $H$, the function resulted is continuous.We note that:\\
$H_{1,1}(u,0)=h(u)$ and $H_{1,2}(u,0)=H_{1,3}(u,0)=H_{1,4}(u,0)=u$,
$H_{1,1}(u {\rm{ \, such\, that }} \,||u||_{E,\beta}=R_n,t)=H_{1,2}(u {\rm{ \, such\, that }} \,||u||_{E,\beta}=R_n,t)=u$\\
The continuity of $h$ implies that $H(u,t)$ is continuous at $t=0$.

Now to show the continuity for $H$ for $0 \leq t\leq 1$ we must show that the \textit{boundary conditions} are preserved, i.e:
\be
H_{1,1}(u {\rm{ \, such\, that }} \,||u||_{E,\beta}=R_n,t)=H_{1,2}(u {\rm{ \, such\, that }} \,||u||_{E,\beta}=2R_n,t)=u
\ee
\be
H_{1,2}(u {\rm{ \, such\, that }} \,||u||_{E,\beta}=2R_n,t)=H_{1,3}(u {\rm{ \, such\, that }} \,||u||_{E,\beta}=2R_n,t)=(1-t)u
\ee
\be
H_{1,3}(u {\rm{ \, such\, that }} \,||u||_{E,\beta}=3R_n,t)=H_{1,4}(u {\rm{ \, such\, that }} \,||u||_{E,\beta}=3R_n,t)=u+te_{n+1},
\ee
 hence we can conclude that $H$ is continuous for $0\leq t\leq 1$. 

To show that $J_\beta(H(u,t))$ is bounded independently of $\beta$ note that $H_{1,1}(u,t)=h(u)$ and that by hypothesis 
$J_\beta(h(u))\leq b_n^m+\frac{\delta}{2}$. Also $H_{1,2}(u,t)\in E^{+n},H_{1,3}(u,t)\in E^{+(n+1)},H_{1,4}(u,t)\in E^{+(n+1)}$ and by lemma \ref{D_n^m} and \ref{R_n} we conclude
\be
J_\beta(H(u,t))\leq c(n+1).
\ee \\
Now if $1\leq t\leq 2:$ We do not have any a priori estimates \textit{ independently of} $m$ on the dimension of the subspace in which $h(u)$ lies. We deform $h(u)$ to the $0$-map thereby ensuring that the range of $H(u,2)$ lies in $E^{+(n+1)}$ at $t=2$, where by lemmas \ref{R_n}, \ref{D_n^m}, $J_\beta$ is bounded independently of $\beta$:\\
Define $h_1(u)=h(u)$ for $||u||_{E,\beta}\leq R_n$ and $h_1(u)=H_{1,2}(u,1)$ for $R_n\leq ||u||_{E,\beta}\leq 2R_n $:
\be
h_1=\left\{ \begin{array}{ll}
h_1(u)=H_{1,1}(u)=h(u) & 0\leq ||u||_{E,\beta}\leq R_n     \\
h_1(u)=H_{1,2}(u,1)=-u+\frac{2R_n}{||u||_{E,\beta}}u & R_n\leq ||u||_{E,\beta}\leq 2R_n \\
\end{array}\right.
\label{h1}\ee

$h_1$ this constructed is continuous as
\be
H_{1,1}(u {\rm{ \, such\, that }} \,||u||_{E,\beta}=R_n,1)=H_{1,2}(u {\rm{ \, such\, that }} \,||u||_{E,\beta}=R_n,1)=u.
\ee Now,
\be
h_1:B(2R_n,E^{+n})\rightarrow E^{+m}\oplus E^{-m}\oplus N^m
\ee
$H_1(u,t)=h_1((2-t)u)$ is continuous: the continuous map $(2-t)u$
\be
(2-t)u=\left\{ \begin{array}{l}
(2-t)u:B(2R_n,E^{+n})\times [1,2]\rightarrow B(2R_n,E^{+n})      \\
(u,t)\rightarrow (2-t)u  \\
\end{array}\right.
\ee
is composed with $h_1$
so the composition
 \be
  h_1 \ o \ (2-t)u=\left\{ \begin{array}{l}
 E^{+n}\times\R\rightarrow E^{+m}\oplus E^{-m}\oplus N^m     \\
 (u,t)\rightarrow  h_1((2-t)u)  \\
 \end{array}\right.
 \ee
 is continuous. 
 \be
 h_1 \ o \ (2-t)u=\left\{ \begin{array}{l}
 E^{+n}\times\R\rightarrow E^{+m}\oplus E^{-m}\oplus N^m      \\
 (u,t)\rightarrow h_1((2-t)u)\\
 \end{array}\right.
 \ee
 
\begin{displaymath}
\left\{ \begin{array}{ll}
0\leq ||u||_{E,\beta}\leq R_n     & H(u,t)=H_{2,1}(u,t)=h_1((2-t)u)  \\
0\leq ||u||_{E,\beta}\leq 2R_n     & H(u,t)=H_{2,2}(u,t)=h_1((2-t)u)  \\
2R_n\leq ||u||_{E,\beta}\leq 3R_n  & H(u,t)=H_{2,3}(u,t)=(3u-\frac{6R_nu}{||u||_{E,\beta}})+t(\frac{||u||}{R_n}-2)e_{n+1} \\
3R_n\leq ||u||_{E,\beta}\leq R_{n+1} & H(u,t)=H_{2,4}(u,t)=u+te_{n+1} \\
\end{array}\right.
\end{displaymath}
$H_{2,1}(u {\rm{ \, such\, that }} \,||u||_{E,\beta}=R_n,t)=H_{2,2}(u {\rm{ \, such\, that }} \,||u||_{E,\beta}=R_n,t)=u$\\
$H_{2,2}(u {\rm{ \, such\, that }} \,||u||_{E,\beta}=2R_n,t)=H_{2,3}(u {\rm{ \, such\, that }} \,||u||_{E,\beta}=2R_n,t)=0$\\
$H_{2,3}(u {\rm{ \, such\, that }} \,||u||_{E,\beta}=3R_n,t)=H_{2,4}(u {\rm{ \, such\, that }} \,||u||_{E,\beta}=3R_n,t)=u+te_{n+1}$\\
$H_{1,1}(u,1)=H_{2,1}(u,1)=u$ for $||u||_{E,\beta}\leq R_n$\\
$H_{1,2}(u,1)=H_{2,2}(u,1)=h_1(u)$ for $R_n\leq ||u||_{E,\beta}\leq 2R_n$ by definition of $h_1$.\\
$H_{1,3}(u,1)=H_{2,3}(u,1)=(3u-\frac{6R_nu}{||u||_{E,\beta}})+(\frac{||u||}{R_n}-2)e_{n+1}$ for $2R_n\leq ||u||_{E,\beta}\leq 3R_n$\\
$H_{1,4}(u,1)=H_{2,4}(u,1)=u+1e_{n+1}$ for $3R_n\leq ||u||_{E,\beta}\leq R_{n+1}$\\
To show that $J_\beta(H(u,t))$ is bounded independently of $\beta,m$ we first consider \\
Case 1: $0\leq ||u||_{E,\beta}\leq 2R_n$ then as $1\leq t\leq 2$ we also have $||(2-t)u||_{E,\beta}\leq 2R_n$ we have now two subcases:\\
a)If $||(2-t)u||_{E,\beta}\leq R_n$, then $h_1((2-t)u)=h((2-t)u)$ and by (\ref{minim}) $J_\beta(h(2-t)u)\leq b_n^m+\delta $,hence $J_\beta(H(u,t))=J_\beta(h(2-t)u)$ is bounded independently of $\beta,m$.\\
b)If $R_n\leq ||(2-t)u||_{E,\beta}\leq 2R_n$, then $h_1((2-t)u)\in E^{+n}$ and by lemma \ref{R_n} and \ref{D_n^m} $J_\beta(h_1(2-t)u)$ is bounded by a constant depending on $n$ and independently of $\beta,m$.\\
Case 2:If $2R_n\leq ||u||_{E,\beta}\leq R_{n+1}$ then $H(u,t)\in E^{+(n+1)}$ and by lemma \ref{D_n^m} and \ref{R_n}  $J_\beta(H(u,t))\leq C(n+1)$, independently of $\beta,m$\\

Now when $2\leq t\leq 3$.
To define $H$ we will first define a function $G$, then $H$ will be a combination of the identity and $G$:
\be
H(u,t)=(3-t)G(u,t)+(t-2)(u+te_{n+1}),\,\, 2\leq t\leq 3
\ee
\begin{displaymath}
\left\{ \begin{array}{ll}
0\leq ||u||_{E,\beta}\leq 2R_n     & G(u,t)=G_1(u,t)\equiv 0  \,{\rm{for}} \, 2\leq t\leq 3 \\
2R_n\leq ||u||_{E,\beta}\leq 3R_n  & G(u,t)=G_2(u,t)=(3u-\frac{6R_nu}{||u||_{E,\beta}})+t(\frac{||u||}{R_n}-2)e_{n+1}\, {\rm{for}}\,  2\leq t\leq 3  \\
3R_n\leq ||u||_{E,\beta}\leq R_{n+1} & G(u,t)=G_3(u,t)=u+te_{n+1} \, {\rm{for}} \, 2\leq t\leq 3\\
\end{array}\right.
\end{displaymath}
$H_{2,1}(u,2)=h_1(0)=0=G_1(u,2)$ for $||u||_{E,\beta}\leq R_n$\\
$H_{2,2}(u,2)=h_1(0)=0=G_1(u,2)$ for $R_n \leq ||u||_{E,\beta}\leq 2R_n$\\
$H_{2,3}(u,2)=G_2(u,2)=(3u-\frac{6R_nu}{||u||_{E,\beta}})+2(\frac{||u||}{R_n}-2)e_{n+1}$ for $2R_n\leq ||u||_{E,\beta}\leq 3R_n$\\
$H_{2,4}(u,2)=G_3(u,2)=u+2e_{n+1}$ for $2R_n\leq ||u||_{E,\beta}\leq 3R_n$.\\
The range of $G(u,t)$ and of the identity map $Id(u,t)=u+te_{n+1}$ both belong to $E^{+(n+1)}$ so the range of $H(u,t)$ is also included  $E^{+(n+1)}$.
We can now extend for all other values of $t$, by 
\be H(u,t)=u+te_{n+1} \in E^{+(n+1)},\,\ 3\leq t\leq R_{n+1}.
\ee
 Such an $H\in\Lambda_n^m(\delta)$ with $J_\beta(H(u,t))$ bounded independently of $\beta,m$ which concludes the proof.

 At this stage we know that by lemma 1.57 in \cite{R82}$, c_n^m(\delta)$ is a critical value if $c_n^m>b_n^m$. Now to show that there is a subsequence $n_q$ such that this is the case we employ the comparison functional $K$ introduced by Tanaka, in lemma 2.2 in \cite{Tanaka88}:
\[
K(w^+)=\frac{1}{2}||w^+||_{E}-\frac{a_0(s)}{s+1}||w^+||_{L^{s+1}}^{s+1},
\]
 where $a_0(s)$ is a positive constant, which satisfies the Palais-Smale condition. The functional $K$ also satisfies the comparison property :
\[
J_{\beta}(w^+)\geq K(w^+)-a_1(f,s)
\]
for any $w^+\in E^+$, $a_1(f,s)$ is a positive constant independent of $\beta,m$. This is a consequence if an application of Young's inequality. We define the minimax sets:
\[
A_n^m=\{ \sigma\in C(S^{m-n},E^{+m}), \sigma (-x)=\sigma (x)\}
\]
where $S^{m-n}\subset E^{+m}$ is the unit sphere in $\R^{m-n+1}$, whose basis consists of eigenvectors $\{e_n,...,e_m\}$. $x\in S^{m-n}$ if and only if
\be
x=\sum_{i=n}^mx_ie_i \,\ {\rm{and}} \,\ \sum_{i=n}^mx_i^2=1
\ee
and the minimax values
\[
\beta_n^m=\sup_{\sigma\in A_n^m}\min_{x\in S^{m-n}}K(\sigma(x))
\]
Properties of the minimax numbers $\beta_n^m$ from \cite{Tanaka88}: There exists sequences $\nu(n),\widetilde{\nu(n)}$
\be
\nu(n)\leq \beta_n^m\leq\widetilde{\nu(n)}
\label{nu}
\ee
such that $\nu(n),\widetilde{\nu(n)}\rightarrow \infty$ as $n\rightarrow\infty$(independently of $m$).

 Borsuk-Ulam type theorem:

\begin{lem}\cite{Tanaka88}Let $a,b\in\N$. Suppose that $h\in C(S^a,\R^{a+b})$, and $g\in C(\R^b,\R^{a+b})$ are continuous mappings such that
\label{Borsuk-Ulam-Tanaka}
\be
h(x)=h(-x) \,\, \rm{for \,\ all } \,\, x\in S^a
\ee
\be
g(-y)=-g(y) \,\, \rm{for \,\, all } \,\, y\in\R^b
\ee
and there is a $r_0$ such that $g(y)=y$ for all $r\geq r_0$.
Then $h(S^a)\cap g(\R^b)\neq \emptyset$
\end{lem}
\begin{lem}\cite{Tanaka88}Let $\gamma\in \Gamma_n^m$ and $\sigma\in A_n^m$, then
\be
[\gamma(D_n^m)\cup\{ u\in E^{+n}\oplus E^{-m}\oplus N^{-m}, ||u||_{\beta,E}\geq R_n \}]\cap \sigma(S^{m-n})\neq\emptyset
\ee
\end{lem}
Proof:
Apply the lemma above with $a=m-n$ and $b=dimension(E^{+n}\oplus E^{-m}\oplus N^{-m})$. Then extend $\gamma$ to all of $ E^{+n}\oplus E^{-m}\oplus N^{-m}$
by extending it by the identity map on $\partial D_n^m$ and view $\sigma(S^{m-n})$ as embedded in $ E^{+m}\oplus E^{-m}\oplus N^{-m}$, then apply the preceding lemma \ref{Borsuk-Ulam-Tanaka}.\\
\begin{lem}$\forall n\in\N$,
\be
b_n^m\geq \beta_n^m -a_1
\ee
where $a_1$ is independent of $n,m,\beta$.
\label{lowerb}
\end{lem}
\begin{lem}(Proposition 4.1\cite{Tanaka88})Suppose that $\beta_n^m<\beta_n^{m+1}$, $m>n+1$, then there exists a $u_n^m\in E^{+m}$ such that
\be
K(u_n^m)\leq \beta_n^m
\label{beta}
\ee
\be
K^\prime\mid_{E^{+m}}(u_n^m)=0
\ee
\be
index K^{\prime\prime}\mid_{E^{+m}}(u_n^m)\geq n
\label{morse}
\ee
\end{lem}
\begin{lem}(Proposition 5.1\cite{Tanaka88}) For any $\varepsilon >0$, there is a constant $C_\varepsilon >0$, such that for $u\in E^+$
\be
index K^{\prime\prime}(u)\leq C_\varepsilon ||u||_{L^{(s-1)(1+\varepsilon)}}^{(s-1)(1+\varepsilon)}
\label{index}
\ee
\end{lem}
\begin{theo} There is a subsequence $n_q$ and $c$ independent of $\beta,m,n$ such that
\be
b_{n_q}^m > n_q^{\frac{s+1}{s}}
\ee
\end{theo}
Proof:\\
The inequality (\ref{nu}) implies that there is a subsequence $n_q$ such that \\
$\beta_{n_q+1}^m>\beta_{n_q}^m$.
\begin{eqnarray}\beta_{n_q}
& \geq & K(u_{n_q}^m)-\frac{1}{2}K^\prime(u_{n_q}^m)u_{n_q}^m \nonumber\\
& \geq & (\frac{1}{2}-\frac{1}{s+1})a_0(s)||u_{n_q}^m||_{s+1}^{s+1}.
\end{eqnarray}
Then for $\varepsilon >0$ small enough
\begin{eqnarray}||u_{n_q}^m||^{s+1}_{s+1}
& \geq & c||u_{n_q}^m||_{(s-1)(1+\varepsilon)}^{s+1}\nonumber\\
& \geq & c[||u_{n_q}^m||_{(s-1)(1+\epsilon)}^{(s-1)(1-\epsilon)}]^{\frac{s+1}{(s-1)(1+\epsilon)}}\\
& \geq & c_\varepsilon {n_q}^{\frac{s+1}{(s-1)(1+\varepsilon)}}
\end{eqnarray}
by combining (\ref{index}) and (\ref{morse}). Now recalling lemma \ref{lowerb} and that for $\varepsilon$ small enough, $\frac{s+1}{(s-1)(1+\varepsilon)}>\frac{s+1}{s}$ the lemma follows.\\
To conclude we recall lemma 1.64 in \cite{R82} which in our case implies that, for $m$ large enough, independently of $\beta$, if $c_n^m=b_n^m$ for all $n\geq n_1$ then $b_n\leq c n^{\frac{s+1}{s}}$. Then by lemma 1.57 in \cite{R82}, $c_{n_q}^m(\delta)$ is a critical value of $I_\beta$ in $E^{+m}\oplus E^{-m}\oplus N^{m}$.\\
\section{Regularity}
\begin{theo} {\rm{Let}} $f$ be $C^2$, {\rm{for}} $n$ {\rm{large enough there is a classical solution}} $u=v+w$ {\rm{of the modified problem}} (\ref{mpapier}) .
\end{theo}
Proof:\\
In this proof the constants may dependent on $\beta$ and $f$ but are independent of $m$.
The proof of this theorem here is slightly simpler from the one in \cite{R84} as we take advantage of the polynomial growth of the nonlinear term and employ Galerkin approximation.\\
Let $u_{n_q}^m=w^m+v^m\in E^{+m}\oplus E^{-m}\oplus N^m$ a distributional solution corresponding to the critical value $c_{n_q}^m(\delta)$, bounded independently of $\beta,m$, and any $\phi\in E^{+m}\oplus E^{-m}\oplus N^m$:
\be
I^\prime(u_{n_q}^m)\phi=0
\label{Galerkin1}\ee
now taking $\phi=v_{tt}^m\in N^m$ we have
\[
(\beta v_{tt}^m,v_{tt}^m)_{L^2}=(|u_{n_q}^m|^{s-1}u_{n_q}^m+f,v_{tt}^m)_{L^2}
\]
\[
\beta||v_{tt}^m||_{L^2}^2\leq |||u_{n_q}^m|^s||_{L^2}||v_{tt}^m||_{L^2}+||f||_{L^2}||v_{tt}^m||_{L^2}
\]
\[
\beta||v_{tt}^m||_{L^2}\leq c||v_{tt}^m||_{L^2}
\]
now by the argument in the proof of the Palais-Smale property we also have
\be
||w^m||_E<c(n_q),\,\ \beta||v_t^m||_{L^2}<c(n_q)
\ee
hence
\[
||v_{tt}^m||_{L^2}\leq c(\beta,f)
\]
we now have
\[
w_{tt}^m-w_{xx}^m=\beta v_{tt}^m+|u^m_n|^{s-1}u^m+f^m(x,t)\in L^2
\]
hence
$w^m\in H^1\cap C^{1}$
by \cite{R67} and \cite{BCN}. This now implies $w^m\in H^2$, $w^m\rightarrow w(\beta)$ pointwise and $w(\beta)\in H^1\cap C^{1}$. Then if $\phi=v_{tttt}^m$ then
\[
(\beta v_{tt}^m,v_{tttt}^m)_{L^2}=(|u_{n_q}^m|^{s-1}u_{n_q}^m+f,v_{tttt}^m)_{L^2}
\]
\[
(\beta v_{ttt}^m,v_{ttt}^m)_{L^2}=([|u_{n_q}^m|^{s-1}u_{n_q}^m+f]_t,v_{ttt}^m)_{L^2}
\]
and we deduce $||v_{ttt}^m||_{L^2}\leq c(\beta,f)$ hence $v^m_{ttt}\rightarrow v_{tt}(\beta)\in C^0$ hence $v(\beta)$ is $C^2$ and $w(\beta)$ is $C^1$ by applying \cite{BCN} to (\ref{mpapier}) . We now have
\be
u_{n_q}^m\rightarrow u(\beta)\in C^1 \,\, {\rm{as}}\,\ m\rightarrow\infty
\label{pointwise}
\ee
and since (\ref{Galerkin1}) holds for any $\phi\in E^{+m}\oplus E^{-m}\oplus N^m$ we can deduce
\be
I^\prime(u(\beta))\phi=0 \,\,\, \forall \phi\in E\oplus N,
\ee
and $u(\beta)$ is a weak solution of (\ref{mpapier}).
Now for any $\phi\in C^{\infty}\cap L^2(S^1)$ we have
\begin{eqnarray} I^\prime(u(\beta))[\phi(x+t)-\phi(x-t)]
& = & \int_Q [-\beta(p^{\prime\prime}(x+t)-p^{\prime\prime}(-x+t)+|u(\beta)|^{s-1}u(\beta))+f(x,t)]\nonumber\\
&   & [\phi(x+t)-\phi(-x+t)]dxdt\nonumber
\end{eqnarray}
Denoting $\psi(x,t):=[-\beta(p^{\prime\prime}(x+t)+|u(\beta)|^{s-1}u(x,t)+f(x,t)]$ and noting that the functions $\psi,\phi$ are periodic we deduce as in
\cite{R78} that
\[
\int_0^{2\pi}\int_0^\pi\psi(x,t)\phi(x+t)dxdt =\int_0^\pi\int_0^{2\pi}\psi(r,r-x)\phi(r)dxdr
\]
and
\[
\int_0^\pi\int_0^{2\pi}\psi(x,t)\phi(-x+t)dxdt=\int_0^\pi\int_0^{2\pi}\psi(x,r+x)\phi(r)dxdr
\]
for all $\phi\in C^{\infty}\cap L^2(S^1)$ hence
\[
\int_0^\pi\psi(x,r+x)-\psi(x,r-x)dxdr=0
\]
and we have
\be
2\pi \beta p^{\prime\prime}(r)=\int_0^{\pi}(|u(\beta)|^{s-1}u(\beta)(x,r-x)-|u(\beta)|^{s-1}u(\beta)(x,r+x))+f(x,r-x)-f(x,r+x)dx
\label{kernelofBox}
\ee
so $p$ is $C^3$ since $u(\beta)\in C^1$. Since RHS of (\ref{mpapier}) is $C^1$
then by \cite{BCN} $w\in C^2$ and $u(\beta)$ is a classical solution of (\ref{mpapier}).
\begin{lem}There is a constant $c$ independent of $\beta,m$ such that
\be
||w(\beta)||_{C^0}\leq c
\ee
\end{lem}
Proof:\\
By (\ref{us+1}), the bound on $c^m_{n_q}(\delta)$ independent of $\beta,m$ and (\ref{pointwise}), we deduce that
 \be
 \int|u(\beta)|^{s+1}dxdt\leq c(n_q) \, \rm{{independent \, of}} \, \beta.
 \ee
Then by (\ref{kernelofBox})
$||\beta v_{tt}||_{L^1}$ is bounded independently of $\beta$, hence by Lovicarova's formula \cite{Lovicarova} we conclude that there is a constant $c$
\be
||w(\beta)||_{C^0}\leq c(n_q)
\ee
which is independent of $\beta$.
\begin{lem}
There is a constant $c(n_q)$, independent of $\beta$ such that
\be
||v(\beta)||_{C^0}\leq c(n_q).
\ee
\end{lem}
Proof:\\
$\forall \phi\in N$,
\[ \int_0^{\pi}\int_0^{2\pi}(-\beta v_{tt}(\beta)+(g(u(\beta))+f(x,t))\phi dxdt =0\]
\be \int_0^{\pi}\int_0^{2\pi}\beta v_t(\beta)\phi_t+(g(v(\beta)+w(\beta))-g(w))\phi dxdt=-\int_0^{\pi}\int_0^\pi (f(x,t))+g(w))\phi dxdt
\label{vinfinity}
\ee
Define $q$:
\begin{equation}
\left\{ \begin{array}{ll}   q(s)=0 ,\,\ \rm{if} \,\ |s|\leq M . & \\
q(s)=s+M \,\, \rm{if} \,\, s\geq M \,\, \rm{and} \,\, q(s)=s-M \,\, \rm{if} \,\, s\leq M. &
\end{array} \right.
\end{equation}
Now define the function $\psi_K(z)$:
\begin{equation}
\left\{ \begin{array}{ll}   \psi_K(z)=\max_{|\xi|\leq M_5}f_K(z+\xi)-f_K(\xi) \,\, \rm{if} \,\, z>0. & \\
\psi_K(z)=-\min_{|\xi|\leq M_5}(f_K(\xi)-f_K(z+\xi))\,\, \rm{if} \,\, z<0 &
\end{array} \right.
\end{equation}
$\psi_K$ is monotonically increasing and $\lim_{z\rightarrow \pm\infty}\psi_K(z)=\pm\infty $. For $z\geq 0$, $\mu(z)=\min(\psi(z),\psi(-z))$. Define
\[ T_\delta=\{(x,t)\in[0,\pi]\times[0,2\pi] \,\,|v(\beta)|\geq \delta \}.\]
By taking the test function $\phi=q(v^+)-q(v^-)=v^+-v^-$ and noting that $g$ is strictly increasing we have the estimate following lemma 3.7 in \cite{R78}:
\be
\int_{T_\delta}(g(v+w)-g(v))(q^+-q^-)dxdt\geq \frac{M-\delta}{||v||_{C^0}}\mu(\delta)\int_{T_\delta}(|q^+|+|q^-|)dxdt
\ee
hence:
\be
(||g(w)||_{C^0}+||f||_{C^0})\int_T|q^+|+|q^-|dxdt\geq \frac{M-\delta}{||v||_{C^0}}\mu(\delta)\int_{T_\delta}(|q^+|+|q^-|)dxdt.
\ee
Denoting $\max(||v^+||_{C^0},||v^-||_{C^0})=||v^{\pm}||_{C^0}$ we have
\be
\mu(\frac{1}{2}||v^{\pm}||_{C^0})\leq 4(||f||_{C^0}+||g(w)||_{C^0})
\ee
and we can conclude that there is a constant $c$ independent of $\beta$ such that
\be
||v(\beta)||_{C^0}\leq c.
\ee
\begin{lem}The family $v(\beta)$ is equicontinuous.
\end{lem}
Proof:
$u=v+w$. Define $\widehat{v}(x,t)=v(x,t+h)$,$\widehat{w}(x,t)=w(x,t+h)$ and $\widehat{u}=\widehat{v}+\widehat{w}$,$\widehat{f}=f(x,t+h)$,$U=V+W$, where $V=\widehat{v}-v$,$W=\widehat{w}-w$, $q(V^+)=Q^+$,$q(V^-)=Q^-$
\be
\int_T\beta V_t\phi_tdxdt +\int_T g(\widehat{v}+w)-g(u)dxdt=-\int_T g(\widehat{u})-g(\widehat{v}+w)+\widehat{f}-fdxdt
\ee
For $\phi=q(V^+)-q(V^-)$ and $V^+=\widehat{v^+}-v^+$, we have
\be
\int_T [g(V+u)-g(u)+\widehat{f}-f][Q^+-Q^-]dxdt\leq (||f(\widehat{u})-f(\widehat{v}+w)||_{C^0}+||\widehat{f}-f||_{C^0})\int_T(|Q^+|+|Q^-|)dxdt
\ee
and
\be
\int_T [g(V+u)-g(u)][Q^+-Q^-]dxdt \geq \frac{\mu(\delta)(M-\delta)}{||V||_{C^0}}\int_T[|Q^+|+|Q^-|]dxdt.
\ee
Since $w(\beta)\in C^1$ and $f\in C^1$ we deduce
\be
||f(\widehat{u})-f(\widehat{v}+w)||_{C^0}+||\widehat{f}-f||_{C^0})\leq c|h|
\ee
where $c$ is independent of $\beta$, thus
\be
\mu(\frac{1}{2}||V^{\pm}||_{C^0})\leq c|h|
\ee
and the modulus of continuity of $v(\beta)$ is independent of $\beta$.
\begin{theo}The problem (\ref{papier}),(\ref{bdry}) has an infinite number of weak solutions $u=w+v$ where $w\in C^1$ and $v\in C^0$.
\end{theo}
Proof:\\
$||\beta v_{tt}||_{L^1}\rightarrow 0$ as $\beta\rightarrow 0$: Recalling the interpolation inequalities \cite{R78},\cite{Nirenberg} and (\ref{kernelofBox}):
\be
\beta||v_{tt}||_{L^1}\leq\beta||v_{tt}||_{C^0}^{\frac{1}{2}}||v(\beta)||_{C^0}^{\frac{1}{2}}\rightarrow 0
\ee
and Lovicarova fundamental solution in \cite{Lovicarova} implies that $w\in C^1$.\\
Case 1:\\
If $\exists \overline{r}$ such that $u(x,\overline{r}-x)=\alpha$ for $\forall x\in[0,\pi]$ then the boundary conditions imply $\alpha=0$ and $p(\overline{r}-2x)=p(\overline{r})+w(x,\overline{r}-x)$, thus
\be
||v||_{C^1}\leq ||w||_{C^1}.
\ee
Case 2:\\ There is no $\overline{r}$ such that $u(x,\overline{r}-x)=0$, then there is $\gamma>0$ such that $\int_0^\pi s|u|^{s-1}(x,r-x)dx>\gamma$,
$\forall r\in [0,2\pi]$. Now since $u(\beta)\rightarrow $ as $\beta\rightarrow 0$ we have
\be
\int_0^\pi s|u|^{s-1}(\beta)(x,r-x)dx>\frac{\gamma}{2}
\ee
Differentiating (\ref{kernelofBox}) with refer to $r$ and using the boundary conditions for $u$ as in \cite{R78} we obtain:
\begin{eqnarray}-\pi\beta p^{\prime\prime\prime}(r)+a(r)p^\prime(r)
& = & \int_0^\pi s|u|^{s-1}(x,r-x)[-\frac{1}{2}w_x(x,r-x)-w_r(x,r-x)]+\nonumber\\
&    & s|u|^{s-1}(x,r+x)[-\frac{1}{2}w_x(x,r+x)+w_r(x,r+x)]+\nonumber\\
&     & f_r(x,r+x)-f_r(x,r-x)dx,
\end{eqnarray}
where $a(r)=\int_0^\pi s|u|^{s-1}(\beta)(x,r-x)+s|u|^{s-1}(\beta)(x,r+x)dx$.
Now by writing $\phi(r)=p^\prime(r)$ we have:
\be
-\pi \beta\phi^{\prime\prime}(r)+a(r)\phi(r)=h(r)
\label{v0}
\ee
where $h\in C^0(S^1)$ and since $f\in C^1$ we deduce as in \cite{R78} that $\lim_{\beta\rightarrow 0}\phi(\beta)$ exists and is in $H^1(S^1)$. Denoting this limit by $\phi(0)$ we deduce that $v\in C^1$. This implies $w\in C^2$ and $h\in C^1$, as $f\in C^2$. Now (\ref{v0}) is valid a.e at $\beta=0$ which implies $\phi\in C^1$ and $u\in C^2$ is a classical solution of (\ref{papier}),(\ref{bdry}).
\bibliographystyle{plain}

\end{document}